\documentclass{amsart}

\usepackage{amssymb,amsfonts,latexsym}

\numberwithin{equation}{section}

\newtheorem{theorem}{Theorem}[section]

\newtheorem{lemma}[theorem]{Lemma}

\theoremstyle{definition}

\newtheorem *{Theorem A}{Theorem A}
\newtheorem *{Theorem B}{Theorem B}
\newtheorem *{Corollary C}{Corollary C}

\newcommand{\btu}{\bar{\tau}}

\newcommand{\tut}{\eta_t}

\newcommand{\Qn}{Q_{2^n}}
\newcommand{\Dn}{D_{2^n}}
\newcommand{\Cn}{C_{2^{n-1}}}

\newcommand{\la}{\langle\,}

\newcommand{\ra}{\,\rangle}

\newcommand{\ut}{\text{Aut}}

\newcommand{\ben}{\begin{enumerate}}

\newcommand{\een}{\end{enumerate}}

\newcommand{\bx}{{\bar x}}

\newcommand{\C}{{\mathbb C}}

\begin{document}

\title
[]
{ Separable deformations of the generalized quaternion group
algebras}

\author{Yuval Ginosar}

\address{Department of Mathematics, University of Haifa, Haifa 31905, Israel}
\email{ginosar@math.haifa.ac.il}

\begin{abstract}
The group algebras $kQ_{2^n}$ of the generalized quaternion groups $\Qn$
over fields $k$ which contain $\mathbb{F}_{2^{n-2}}$, are deformed to separable $k((t))$-algebras $[kQ_{2^n}]_t$.
The dimensions of the simple components of $\overline{k((t))}\otimes_{k((t))}[kQ_{2^n}]_t$
over the algebraic closure $\overline{k((t))}$, and those of $\C Q_{2^n}$ over $\C$ are the same,
yielding strong solutions of the Donald-Flanigan conjecture for the generalized quaternion groups.
\end{abstract}
\maketitle

The Donald-Flanigan (DF) conjecture \cite{DF} says that any group algebra $kG$ of a finite group $G$ over a field $k$ admits a separable deformation.
It was verified in \cite{E,ES,GG,GS96,GS97,GGS,M,PS,S} for certain families of finite groups.
In \cite{BG} a separable deformation was constructed for the quaternion group $G=Q_8$, turning the generalized quaternion group $Q_{16}$ to a current minimal unsolved case.
In this note, we extend the strategy of \cite{BG} in order to deal, as promised therein, with the family of generalized quaternion groups $\Qn$.
As a by-product we establish separable deformations for the family of dihedral 2-groups $\Dn$
(K. Erdmann and M. Schaps have already found separable deformations for this family in \cite{ES}).
Certainly, for both families of 2-groups, the interesting case is where the field $k$ is of characteristic 2.
For our considerations, $k$ is further assumed to contain the Galois field $\mathbb{F}_{2^{n-2}}$ of $2^{n-2}$ elements.
The solutions to the DF conjecture for $\Qn$ and $\Dn$ established here are so called strong.

\section{Background}\label{bground}
Let $A$ be a $k$-algebra (all the algebras throughout are associative), let $k[[t]]$ be the ring of formal power series over $k$, and let $k((t))$ be its field of fractions.
Suppose that the free $k[[t]]$-module $\Lambda_t:=k[[t]]\otimes_k A$ admits a multiplication such that there is an isomorphism $A\simeq \Lambda_t/\langle t\Lambda_t\rangle$
of $k$-algebras. Then the $k((t))$-algebra
$A_t:= k((t))\otimes_{k[[t]]}\Lambda_t$ is called a {\it deformation} of $A$. The algebra $A$ specializes $\Lambda_t$ at $t=0$.
Even though $t$ is invertible in a deformation $A_t$, we adopt an abuse of language saying that $A$ is a specialization also of $A_t$ at $t=0$.

Consider the extension
\begin{equation}\label{extension}
[\beta]: 1\to \Cn\to \Qn\to C_2\to 1,
\end{equation}
with an action $\eta:C_2\to \ut (\Cn)$ of $C_2=\la \bar\tau\ra  $ on $\Cn=\la \sigma\ra  $ via
\begin{equation}\label{act}
\eta(\btu):\sigma\mapsto \sigma^{2^{n-1}-1}(=\sigma^{-1}),
\end{equation}
and an associated 2-cocycle $\beta:C_2\times C_2\to \Cn$ representing \eqref{extension} which is given by
$$\beta(1,1)=\beta(1,\btu)=\beta(\btu,1)=1,\beta(\btu,\btu)=\sigma^{2^{n-2}}.$$
The group algebra $k\Qn$ can be viewed as follows.
First, the group automorphism \eqref{act} of $\Cn$ can be extended to an algebra automorphism of the subgroup algebra $k\Cn$ taking $\sigma$ to its inverse.
Next, let $k\Cn[y;\eta]$ be the skew polynomial ring (see \cite[\S 1.2]{MR}) over $k\Cn$,
where we keep the notation $\eta$ for the action of the indeterminate $y$ on $k\Cn$ via the above extension of \eqref{act}. Then $k\Qn$
is isomorphic to the quotient $k\Cn[y;\eta]/\la g(y)\ra  $ of this skew polynomial ring by the central polynomial
\begin{equation}\label{qn}
g(y):=y^{2}-\sigma^{2^{n-2}}  \in k\Cn[y;\eta].
\end{equation}
The base algebra $k\Cn$ can itself be identified with a quotient
\begin{eqnarray}\label{ident}
\begin{array}{rcl}
k\Cn&\xrightarrow{\simeq} &k[x]/\la x^{2^{n-1}}-1\ra\\
\sigma&\mapsto&\bx:=x+\la x^{2^{n-1}}-1\ra.
\end{array}
\end{eqnarray}
The above description is good for any field $k$. From now on, the above condition $\mathbb{F}_{2^{n-2}}\subseteq k$ is entailed.

Here is a layout of the paper. In
\S \ref{C4.1} the subgroup algebra $kC_{2^{n-1}}$ is deformed to a
separable algebra $[k{C_{2^{n-1}}}]_t$ which is isomorphic to a direct sum of $2^{n-2}+1$ fields
$\bigoplus_{a\in \mathbb{F}^*_{2^{n-2}}}K_a\oplus k((t))\oplus k((t))$, where $K_a$ are separable field extensions
of $k((t))$ of degree 2.
The next step (\S\ref{C4.2}) is to
construct an automorphism $\tut$ of $[k{C_{2^{n-1}}}]_t$ which agrees
with the action $\eta$ of $C_2$ on $kC_{2^{n-1}}$ when specializing  $t=0$.
This action fixes all the primitive idempotents of $[{kC_{2^{n-1}}}]_t$.
%This action will interchange two of the copies of $k((t))$
%and fix the two other copies.
By that we obtain the skew polynomial ring $[k{C_{2^{n-1}}}]_t[y;\tut]$.
In \S \ref{hbt} we deform the polynomial $g(y)$ \eqref{qn} to
a separable polynomial $g_t(y)$ of degree 2 in $y$, which lies in the center of $[kC_{2^{n-1}}]_t[y;\tut]$.
Factoring out the two-sided ideal generated by $g_t(y)$, we establish a deformation
$$[{kQ_{2^{n}}}]_t:=[{kC_{2^{n-1}}}]_t[y;\tut]/\la g_t(y)\ra  $$ of
$kQ_{2^{n}}={kC_{2^{n-1}}}[y;\eta]/\la g(y)\ra$. The proofs of the claims of \S\ref{defor} are postponed and given in \S\ref{proofs}.
In \S\ref{sep} we show that $[{kQ_{2^{n}}}]_t$ as above is separable. In \S\ref{dihed} we adapt the above strategy to the $2^n$-dihedral group algebra $k\Dn$,
constructing a separable deformation $[{kD_{2^{n}}}]_t$.
Moreover, passing to the algebraic closure $\overline{k((t))}$
we have
$$\overline{k((t))}\otimes_{k((t))} [kQ_{2^{n}}]_t\simeq\overline{k((t))}\otimes_{k((t))} [kD_{2^{n}}]_t\simeq \bigoplus_{i=1}^4\overline{k((t))}\oplus
\bigoplus\limits_{i=1}^{2^{n-2}-1} M_2(\overline{k((t))}).$$
These are {\it strong} solutions to the DF conjecture since
their decompositions to simple components afford the same dimensions as
$$\C \Qn\simeq\C\Dn\simeq\bigoplus\limits_{i=1}^4\C\oplus \bigoplus\limits_{i=1}^{2^{n-2}-1}M_2(\C).$$

\section{The deformation}\label{defor}
\subsection{}\label{C4.1}
We first deform the cyclic group algebra $k\Cn$
%(where $k$ contains $\mathbb{F}_{2^{n-2}}$)
as follows.
Let $$U:=\{1+tk[[t]]\}\subseteq k[[t]]$$ be the group of {\it 1-units} of $k[[t]]$.
For any distinct 1-units $c_1,c_2\in U$ (which are determined more precisely in the next section), define
\begin{equation}\label{pitx}
\pi_t(x):=(x+c_1)(x+c_2)\prod_{a\in \mathbb{F}^*_{2^{n-2}}}(x^2+atx+1)\in k((t))[x].
\end{equation}
Then $\pi_{t=0}(x)=x^{2^{n-1}}+1$, and hence by \eqref{ident}, the algebra
$$[kC_{2^{n-1}}]_t:= k((t))[x]/\la \pi_t(x)\ra $$ is a deformation of $kC_{2^{n-1}}$.
Note that for any element $a\in \mathbb{F}_{2^{n-2}}^*(=\mathbb{F}_{2^{n-2}}\setminus \{0\})$,
the polynomial $x^2+atx+1\in k((t))[x]$ is separable. Furthermore, it does not admit roots mod($t^2$), and is hence irreducible. The polynomial $\pi_t(x)$
is given as a product \eqref{pitx} of distinct irreducible polynomials, each of which is separable, and so is itself separable.
The commutative algebra $[kC_{2^{n-1}}]_t$ is then a sum of fields, corresponding to the irreducible factors of $\pi_t(x)$ as follows.
Denote the primitive idempotents of $[kC_{2^{n-1}}]_t$ by
\begin{equation}\label{idemps}
e_{c_1}\in\left \la\frac{\pi_t(\bx)}{\bx+c_1}\right \rangle,\ \ e_{c_2}\in\left \la\frac{\pi_t(\bx)}{\bx+c_2}\right \rangle,\ \
e_a\in\left \la\frac{\pi_t(\bx)}{\bx^2+at\bx+1}\right \rangle,\ \ \forall a\in \mathbb{F}^*_{2^{n-2}}.
\end{equation}
%$e_{c_1}, e_{c_2}$ and $e_a$ for every ${a\in \mathbb{F}^*_{2^{n-2}}}$
The fields
$$k_1:=[kC_{2^{n-1}}]_te_{c_1},\ \ k_2:=[kC_{2^{n-1}}]_te_{c_2}$$ are one-dimensional over $k((t))$, and
the fields $$K_a:=[kC_{2^{n-1}}]_te_a\simeq k((t))[x]/\la x^2+atx+1\ra$$
are two-dimensional over $k((t))$ for every ${a\in \mathbb{F}^*_{2^{n-2}}}$. Write
\begin{equation}\label{fieldec}
[kC_{2^{n-1}}]_t=
%\bigoplus_{a\in \mathbb{F}^*_{2^{n-2}}}[kC_{2^{n-1}}]_te_a\oplus [kC_{2^{n-1}}]_te_{c_1}\oplus [kC_{2^{n-1}}]_te_{c_2}\simeq
\bigoplus_{a\in \mathbb{F}^*_{2^{n-2}}}K_a\oplus k_1\oplus k_2.
\end{equation}
As customary, denote $\bx:=x+\la \pi_t(x)\ra \in[kC_{2^{n-1}}]_t$.
We record the following claim for a later use, it is proven in \S\ref{caproof}.
\begin{lemma}\label{ca}
For every $a\in \mathbb{F}^*_{2^{n-2}}$ there exists $d_a\in k[[t]]\setminus\{0\}$
such that
$$(\bx^{2^{n-2}}+t^{2^{n-2}-1}\bx) e_a=d_ae_a.$$
\end{lemma}

\subsection{}\label{C4.2}
Our next step is to deform the action $\eta$.
Let
\begin{equation}\label{pt}
p_t(x):=(dx^2+btx+1)\prod_{a\in\mathbb{F}_{2^{n-2}}^*}(x^2+atx+1)\in k((t))[x]
\end{equation}
%(here $a$ runs over {\it all} the elements of $\mathbb{F}_{2^{n-2}}$).
for some $d\in U$ and $b\in k[[t]]$ (both are later to be chosen). Then it is not hard to verify that
%as a product of reciprocal polynomials, $p_t(x)$ is itself reciprocal of degree $2^{n-1}$, and satisfies
\begin{equation}\label{pt0}
p_{t=0}(x)=x^{2^{n-1}}+1.
%=(x+1)^{2^{n-1}}.
\end{equation}
\begin{lemma}\label{distinct1units}
With the notation \eqref{pt} there exist $d\in U$ and $b\in k[[t]]$ such that the polynomial
\begin{equation}\label{ft}
f_t(x):=p_t(x)+x^2+1\in k((t))[x]
\end{equation} admits two distinct roots in $U$.
\end{lemma}

Owing to Lemma \ref{distinct1units}, whose proof can be found in \S\ref{cd}, we choose $b\in k[[t]]$ and $d\in U$ such that $c_1, c_2\in U$ are distinct roots of $f_t(x)$.
These are the 1-units in $\pi_t(x)$ \eqref{pitx}. We record
\begin{equation}\label{ftci}
f_t(c_i)=p_t(c_i)+c_i^2+1=0,\ \ i=1,2.
\end{equation}

Note that since $p_t(0)=1$, then $p_t(x)+1\in xk((t))[x]$ (the ideal generated by $x$ in the polynomial algebra $k((t))[x]$), and hence
\begin{equation}\label{qt}
\tilde{\eta}_t(x):=\frac{p_t(x)+1}{x}\in k((t))[x].
\end{equation}
The polynomial \eqref{qt} determines a $k((t))$-algebra endomorphism
\begin{eqnarray}\label{p8}
\tilde{\tut}:\begin{array}{ccc}
k((t))[x]&\rightarrow &k((t))[x]\\
x&\mapsto & \tilde{\eta}_t(x)
\end{array}.
\end{eqnarray}
\begin{lemma}\label{tut}
With the notation above, $\tilde{\tut}$ induces an automorphism
\begin{eqnarray}\label{eta}
\tut:\begin{array}{ccc}
k((t))[x]/\la \pi_t(x)\ra&\rightarrow &k((t))[x]/\la \pi_t(x)\ra\\
\bx&\mapsto & \overline{\tilde{\eta}_t(x)}
\end{array},
\end{eqnarray}
of order 2, which fixes all the idempotents of $[kC_{2^{n-1}}]_t= k((t))[x]/\la \pi_t(x)\ra$
\begin{equation}\label{ei}
\tut(e_{c_1})=e_{c_1},\ \  \tut(e_{c_2})=e_{c_2},\ \  \tut(e_a)=e_a,\ \ \forall a\in \mathbb{F}^*_{2^{n-2}},
\end{equation}
furthermore, for every $a\in \mathbb{F}^*_{2^{n-2}}$
\begin{equation}\label{txe}
\tut(\bx e_a)=(\bx+at)e_a.
\end{equation}
\end{lemma}
By Lemma \ref{tut}, whose proof can be found in \S\ref{inidpi}, $\tut$ induces automorphisms of order 2
of the fields $K_a\subseteq [kC_{2^{n-1}}]_t$ while fixing the two fields $k_1,k_2\subseteq [kC_{2^{n-1}}]_t$ pointwise.
Furthermore, by the definitions \eqref{qt}, \eqref{p8}, and equation \eqref{pt0} we have
$$\eta_{t=0}:\bx\mapsto\frac{p_{t=0}(\bx)+1}{\bx}=\bx^{2^{n-1}-1}.$$
Consequently, the automorphism
$\eta_{t=0}$ of $[{kC_{2^{n-1}}}]_{t=0}$ agrees with the automorphism $\eta(\btu)$ of $kC_{2^{n-1}}$ (with the identification \eqref{ident}).
The skew polynomial ring
$$[{kC_{2^{n-1}}}]_t[y;\tut]=(k((t))[x]/\la \pi_t(x)\ra) [y;\tut]$$
is therefore a deformation of $kC_{2^{n-1}}[y;\eta]$.
Note that by (\ref{ei}), the
idempotents $e_{c_1},e_{c_2}, e_a (a\in \mathbb{F}^*_{2^{n-2}})$ are central in $[{kC_{2^{n-1}}}]_t[y;\tut]$ and hence
\begin{equation}\label{sumid}
[{kC_{2^{n-1}}}]_t[y;\tut]=\bigoplus\limits_{a\in \mathbb{F}^*_{2^{n-2}}}[{kC_{2^{n-1}}}]_t[y;\tut]e_a\oplus
[{kC_{2^{n-1}}}]_t[y;\tut]e_{c_1}\oplus[{kC_{2^{n-1}}}]_t[y;\tut]e_{c_2}.
\end{equation}

\subsection{}\label{hbt}
We complete the construction of $[kQ_{2^n}]_t$ by deforming the polynomial $g(y)$ (\ref{qn}).
For any $z\in k[[t]]$ let
\begin{equation}\label{qt1}
g_t(y):=y^2+z(e_{c_1}+e_{c_2})y+\bx^{2^{n-2}}+t^{2^{n-2}-1}\bx\in [kC_{2^{n-1}}]_t[y;\tut].
\end{equation}
%Decomposition of (\ref{qt1}) with respect to the idempotents $e_1,e_2,e_3$ yields
%\begin{equation}\label{qt} q_t(y)=(y^2+b)e_1+[y^2+zay+c(c+a)]e_2+[y^2+zay+d(d+a)]e_3.\end{equation}
\begin{lemma}\label{gcentral}
With the notation \eqref{qt1}, $g_t(y)$ is in the center of $[kC_{2^{n-1}}]_t[y;\tut].$
Consequently, $\la g_t(y)\ra=g_t(y)[kC_{2^{n-1}}]_t[y;\tut]$ is a two-sided ideal.
\end{lemma}
The proof of Lemma \ref{gcentral} is given in \S\ref{gcentproof}. For every $z\in k[[t]]$ the element $z(e_{c_1}+e_{c_2})$ lies in $[kC_{2^{n-1}}]_t$.
Choose a non-zero $z\in k[[t]]$ such that $$z(e_{c_1}+e_{c_2})_{t=0}=0$$ ($z$ can be taken as $t^m$ for
sufficiently large $m$). Plugging this choice of $z$ in \eqref{qt1} and identifying $\sigma$ and $\bx$ as in \eqref{ident} we have
\begin{equation}\label{qt0}
g_{t=0}(y)=y^2+\bx^{2^{n-2}}=g(y).
\end{equation}
%where the leading term $y^2$ remains unchanged.
Lemma \ref{gcentral} and equation \eqref{qt0} yield that
\begin{equation}\label{thedeform}
[kQ_{2^{n}}]_t:=[kC_{2^{n-1}}]_t[y;\tut]/\la g_t(y)\ra
\end{equation}
is a deformation of $kQ_{2^{n}}$.

\section{Separability of $[kQ_{2^n}]_t$}\label{sep}

Separability of the deformed algebra $[kQ_{2^n}]_t$ is proved in the same fashion as in \cite{BG} for the case $n=3$.
By \eqref{fieldec}, \eqref{sumid} and \eqref{thedeform},
\begin{equation}\label{e4}
[kQ_{2^{n}}]_t=\bigoplus\limits_{a\in \mathbb{F}^*_{2^{n-2}}}K_a[y;\tut]/\la g_t(y)e_a \ra\oplus k_1[y;\tut]/\la g_t(y)e_{c_1} \ra\oplus
k_2[y;\tut]/\la g_t(y)e_{c_2} \ra.
\end{equation}
We now show that all the direct summands of \eqref{e4} are separable.

Let $a\in \mathbb{F}^*_{2^{n-2}}$. The non-zero element $d_a\in k[[t]]$ provided in Lemma \ref{ca}, as well as orthogonality of the central idempotents
$e_a$, $e_{c_1}$ and $e_{c_2}$ yield
$$g_t(y)e_a=(y^2+\bx^{2^{n-2}}+t^{2^{n-2}-1}\bx) e_a=(y^2+d_a) e_a.$$
We obtain
$$[{kC_{2^{n-1}}}]_t[y;\tut]e_a/\la g_t(y)e_a\ra=K_a[y;\tut]/\la y^2+d_a \ra\simeq K_a*_{\tut}^{\varphi_a}C_2.$$
The rightmost term is the {\it crossed product} of the group $C_2:=\la \bar{\tau}\ra$
acting faithfully on the field $K_a$ via $\tut$
(\ref{txe}), with a twisting determined by the 2-cocycle
$$\begin{array}{c}
\varphi_a:C_2\times C_2\to K_a^*,\\
 \varphi_a(1,1)=\varphi_a(1,\btu)=\varphi_a(\btu,1)=1,\ \, \varphi_a(\btu,\btu)=d_a.
 \end{array}$$
This is a central simple algebra over the subfield of invariants $k((t))$ \cite[Theorem 4.4.1]{H}.
Evidently, this simple algebra is split by $\overline{k((t))}$, i.e.
\begin{equation} \label{e1}
\overline{k((t))}\otimes_{k((t))}K_a[y;\tut]/\la g_t(y)e_a\ra  \simeq M_2(\overline{k((t))}).
\end{equation}

Next, by Lemma \ref{tut}, the action $\tut$ is trivial on both $k_1$ and $k_2$, hence we may regard the skew polynomial rings
$k_1[y;\tut]$ and $k_2[y;\tut]$ as ordinary polynomial rings $k_1[y]$ and $k_2[y]$ respectively.
Equation \eqref{idemps} yields
$(\bx+c_i) e_{c_i}=0$ for $i=1,2,$
in other words
\begin{equation}\label{eq}
\bx e_{c_1}=c_1e_{c_1},\ \ \bx e_{c_2}=c_2e_{c_2}.
\end{equation}
Orthogonality of the idempotents $e_{c_1}$ and $e_{c_2}$, together with \eqref{eq} yields
$$g_t(y)e_{c_i}= (y^2+zy+{c_i}^{2^{n-2}}+t^{2^{n-2}-1}{c_i})e_{c_i}, \ \ i=1,2.$$
We obtain
$$k_1[y;\tut]/\la g_t(y)e_{c_1}\ra\simeq k[y]/\la y^2+zy+{c_1}^{2^{n-2}}+t^{2^{n-2}-1}{c_1}\ra,$$
and
$$k_2[y;\tut]/\la g_t(y)e_{c_2}\ra\simeq k[y]/\la y^2+zy+{c_2}^{2^{n-2}}+t^{2^{n-2}-1}{c_2}\ra.$$
The polynomials $y^2+zy+{c_1}^{2^{n-2}}+t^{2^{n-2}-1}{c_1}$ and $y^2+zy+{c_2}^{2^{n-2}}+t^{2^{n-2}-1}{c_2}$ in $k[y]$ are separable
(since $z$ is non-zero) and split as products of degree-1 polynomials over the algebraic closure $\overline{k((t))}$.
Both $k_1[y;\tut]/\la g_t(y)e_{c_1}\ra$ and
$k_2[y;\tut]/\la g_t(y)e_{c_2}\ra$ are thus separable
$k((t))$-algebras and
\begin{eqnarray} \label{e23}\begin{array}{l}
\overline{k((t))} \otimes_{k((t))}k_1[y;\tut]/\la g_t(y)e_{c_1}\ra\simeq\\
\simeq \overline{k((t))} \otimes_{k((t))} k_2[y;\tut]/\la g_t(y)e_{c_2}\ra
\simeq \overline{k((t))}\oplus \overline{k((t))}.
\end{array}
\end{eqnarray}
By \eqref{e4},\eqref{e1} and \eqref{e23}, a strong solution is established
$$\overline{k((t))}\otimes_{k((t))}[kQ_{2^{n}}]_t \simeq\bigoplus_{i=1}^{{2^{n-2}-1}}M_2(\overline{k((t))})
\oplus \bigoplus_{i=1}^4\overline{k((t))}.$$

\section{Proofs}\label{proofs}
\subsection{Proof of Lemma \ref{ca}}\label{caproof}
Since the idempotent $e_a$ lies in the ideal generated by $\frac{\pi_t(\bx)}{\bx^2+at\bx+1}$ \eqref{idemps}
% for every $a\in \mathbb{F}^*_{2^{n-2}}$, $e_a\in\la.
it is enough to prove that for every $a\in \mathbb{F}^*_{2^{n-2}}$ there exists $d_a\in k[[t]]\setminus\{0\}$ such that
\begin{equation}
x^{2^{n-2}}+t^{2^{n-2}-1}x\equiv d_a\text{ mod}(x^2+atx+1).
\end{equation}\label{indeed}
Indeed, for any non-negative integer $j$
$$(x^2+atx+1)^{2^j}=x^{2^{j+1}}+(atx)^{2^{j}}+1,$$
and so
\begin{equation}\label{mod}
x^{2^{j+1}}\equiv(atx)^{2^{j}}+1\ \ \mbox{mod}(x^2+atx+1),\ \ \forall j\geq 0.
\end{equation}
Putting $j=n-3$ and then $j=n-4$ in \eqref{mod} we obtain for every $a\in \mathbb{F}^*_{2^{n-2}}$
$$x^{2^{n-2}}\equiv (atx)^{2^{n-3}}+1\equiv (at)^{2^{n-3}}((atx)^{2^{n-4}}+1)+1\ \ \mbox{mod}(x^2+atx+1).$$
Proceeding the iteration of (\ref{mod}) yields
$$x^{2^{n-2}}\equiv (at)^{2^{n-3}}((at)^{2^{n-4}}...((at)^2(atx+1)+1)...+1)+1\equiv(at)^{2^{n-2}-1}x+d_a\ \ \mbox{mod}(x^2+atx+1),$$ where $d_a$ is a
sum of certain powers of $at$, and hence does not depend on $x$. We now
make use of the fact that the elements $a\in\mathbb{F}^*_{2^{n-2}}$ satisfy $a^{2^{n-2}-1}=1$. Consequently,
$$x^{2^{n-2}}\equiv t^{2^{n-2}-1}x+d_a\ \ \mbox{mod}(x^2+atx+1),$$  and \eqref{indeed} is obtained.\qed

\subsection{Proof of Lemma \ref{distinct1units}}\label{cd}
For the sake of simplicity we denote
\begin{equation}\label{mu}
\mu(x):=\prod_{a\in \mathbb{F}_{2^{n-2}}^*}(x^2+atx+1)\in k((t))[x].
\end{equation}
Note that
\begin{equation}\label{equ}
u^2+atu+1\in at\cdot U,\ \ \forall u\in U, a\in \mathbb{F}^*_{2^{n-2}}.
\end{equation}
By \eqref{equ}, using the fact that the product of invertible elements in $\mathbb{F}_{2^{n-2}}$ equals 1, we have
\begin{equation}\label{muu}
\mu(u)\in \prod_{a\in \mathbb{F}_{2^{n-2}}^*}(at\cdot U)=t^{2^{n-2}-1}\cdot U,\ \ \forall u\in U,
\end{equation}
which we record for a later use.

We prove the lemma by an example of two distinct 1-units which annihilate the polynomial $f_t(x)$ for certain
\begin{equation}\label{bd}
b\in k[[t]]\ \ \wedge \ \ d\in U.
\end{equation}
Our 1-units are $1+\alpha_1t^m$ and $1+\alpha_2t^m$, where $m\geq 2^{n-2}$ and where $\alpha_1,\alpha_2\in k^*$ are any distinct elements.
With the notation \eqref{pt}, \eqref{ft} and \eqref{mu}, the elements $1+\alpha_1t^m$ and $1+\alpha_2t^m$ are roots of $f_t(x)$ if and only if for $i=1,2$
$$f_t(1+\alpha_it^m)=[d(1+\alpha_it^{m})^2+bt\cdot(1+\alpha_it^m)+1]\cdot\mu(1+\alpha_it^m)+(1+\alpha_it^{m})^2+1=0,$$
that is, if and only if
\begin{equation}\label{system}
d(1+\alpha_i^2t^{2m})+bt\cdot(1+\alpha_it^m)=\frac{\alpha_i^2t^{2m}}{\mu(1+\alpha_it^m)}+1,\ \ i=1,2.
\end{equation}
Consider \eqref{system} as a system of two non-homogeneous linear equations in the variables $b$ and $d$.
We show that the system admits a (unique) solution satisfying \eqref{bd}.
Indeed, solving the system \eqref{system} for $b$, it can be verified using \eqref{muu} that
%$$b(\alpha_1+\alpha_2)t^{m+1}=t^{2m},$$\left (\frac{\alpha_1^2}{\mu(1+\alpha_1t^m)}+\frac{\alpha_2^2}{\mu(1+\alpha_2t^m)}\right )t^{2m}+(\alpha_1^2+\alpha_2^2)t^{2m},$$and so
%\begin{equation}\label{system1}b=\frac{1}{\alpha_1+\alpha_2}\cdot\left (\frac{\alpha_1^2}{\mu(1+\alpha_1t^m)}+\frac{\alpha_2^2}{\mu(1+\alpha_2t^m)}\right )t^{m-1}+(\alpha_1+\alpha_2)t^{m-1}.
%\end{equation}From and \eqref{system1} we deduce that
$$b\in (\alpha_1+\alpha_2)\cdot t^{m-2^{n-2}}\cdot U.$$ Our choice of $m\geq 2^{n-2}$ ascertains that the left condition in \eqref{bd} is fulfilled.
Returning to \eqref{system}, using the fact that $b\in k[[t]]$ was just established, we get %it is easy to extract $d$ and
$$d(1+\alpha_i^2t^{2m})=bt\cdot(1+\alpha_it^m)+\frac{\alpha_i^2t^{2m}}{\mu(1+\alpha_it^m)}+1\in U.$$
Since $1+\alpha_i^2t^{2m}\in U$, we deduce that $d$ satisfies the right condition in \eqref{bd}.
\qed\\

\subsection{Proof of Lemma \ref{tut}}\label{inidpi}
%Since the algebra $k((t))[x]/\la \pi_t(x)\ra  $ is semisimple, then
We show that $\tilde{\eta}_t$ takes the ideal generated by each irreducible factor of $\pi_t(x)$ to itself.
First, by the definition of $\tilde{\eta}_t$ \eqref{qt},
$$\tilde{\eta_t}(x+c_i)=\tilde{\eta_t}(x)+c_i=\frac{p_t(x)+1}{x}+c_i,\ \ i=1,2.$$
Then by \eqref{ftci} it follows that $c_1$ and $c_2$ annihilate the polynomials $\tilde{\eta}_t(x+c_1)$ and $\tilde{\eta}_t(x+c_2)$ respectively. Hence
%This argument, applied also for $(x+d)$, yields
\begin{equation}\label{inidcd}
\tilde{\eta_t}(x+c_1)\in \la x+c_1\ra,\ \ \tilde{\eta_t}(x+c_2)\in \la x+c_2\ra.
\end{equation}
Next, for every $a\in \mathbb{F}^*_{2^{n-2}}$ develop
\begin{equation}\label{x2tildeta}
x^2\cdot\tilde{\eta_t}(x^2+atx+1)=x^2\cdot\left (\frac{p_t^2(x)+1}{x^2}+at\frac{p_t(x)+1}{x}+1\right )=p_t^2(x)+1+at(xp_t(x)+x)+x^2.
\end{equation}
Computing \eqref{x2tildeta} modulo $x^2+atx+1$, bearing in mind that by the definition \eqref{pt}
\begin{equation}\label{pteq0}
p_t(x)\equiv 0\ \ \mbox{mod}(x^2+atx+1),
\end{equation}
we get
$$x^2\tilde{\eta_t}(x^2+atx+1)\equiv 1+atx +x^2\equiv 0\ \ \mbox{mod}(x^2+atx+1).$$
Since the polynomial $x^2$ is prime to $x^2+atx+1$ we obtain
\begin{equation}\label{inida}
\tilde{\eta_t}(x^2+atx+1)\in\la x^2+atx+1\ra,\ \ \forall a\in \mathbb{F}^*_{2^{n-2}}.
\end{equation}
From equations \eqref{inidcd} and \eqref{inida} we conclude that
$\tilde{\eta_t}(\pi_t(x))\in \la \pi_t(x)\ra$ and thus \eqref{eta} is a well-defined $k((t))$-algebra morphism.
Moreover, by \eqref{inidcd} and \eqref{inida}, all the minimal ideals of $[kC_{2^{n-1}}]_t$, namely
$\la\frac{\pi_t(\bx)}{\bx+c_1}\ra, \la\frac{\pi_t(\bx)}{\bx+c_2}\ra,$ and $\la\frac{\pi_t(\bx)}{\bx^2+at\bx+1}\ra_{ a\in \mathbb{F}^*_{2^{n-2}}}$
are stable under $\eta_t$. Hence, every primitive idempotent (see \eqref{idemps}) is either fixed or vanishes under $\eta_t$.
Since these primitive idempotents sum up to 1, and since $\eta_t(1)=1$, it follows that all are fixed by $\eta_t$ proving \eqref{ei} and, in particular, that $\eta_t$ is an automorphism of
$[kC_{2^{n-1}}]_t$.

Finally, for every $a\in \mathbb{F}^*_{2^{n-2}}$
\begin{equation}\label{txe1}
\tut(\bx e_a)=\frac{p_t(\bx)+1}{\bx}\cdot\tut(e_a).
\end{equation}
We apply \eqref{ei}, which has just been verified, together with \eqref{pteq0} and \eqref{txe1} to obtain
$\tut(\bx e_a)=\frac{1}{\bx}\cdot e_a=(\bx+at)e_a$ proving \eqref{txe}.
By equations \eqref{ei} and \eqref{txe} the automorphism $\tut^2$ is an identity on the elements
$$\{e_{c_1},e_{c_2}\}\cup \{e_a\}_{a\in \mathbb{F}^*_{2^{n-2}}}\cup\{\bx e_a\}_{a\in \mathbb{F}^*_{2^{n-2}}}$$
which form a $k((t))$-basis for $k((t))[x]/\la \pi_t(x)\ra$. Hence $\tut$ is of order 2.
\qed

\subsection{Proof of Lemma \ref{gcentral}}\label{gcentproof}
We show that each of the terms of $g_t(y)$, namely the leading term, the free term and the term $z(e_{c_1}+e_{c_2})y$ lie in the center of $[kC_{2^{n-1}}]_t[y;\tut].$
First, the leading term $y^2$ is central since by Lemma \ref{tut} the automorphism $\tut$ is of order 2.
Next, in order to prove that the free (of $y$) term $\bx^{2^{n-2}}+t^{2^{n-2}-1}\bx\in [kC_{2^{n-1}}]_t$ is central, it is enough to show that multiplying it with all the idempotents
of $[kC_{2^{n-1}}]_t$, i.e., $e_{c_1},e_{c_2}$ and $e_a$ (for every $a\in \mathbb{F}_{2^{n-2}}$) yield $\tut$-invariant elements in $[kC_{2^{n-1}}]_t$.
This is clear for $e_{c_1},e_{c_2}$ since by Lemma \ref{tut}, the subspace Span$_{k((t))}\{e_{c_1},e_{c_2}\}=k_1\oplus k_2$ is $\tut$-invariant.
As for the idempotents $e_a$,
Lemma \ref{ca} says that for every $a\in \mathbb{F}^*_{2^{n-2}}$ the projection $(\bx^{2^{n-2}}+t^{2^{n-2}-1}\bx)e_a$ is equal to $d_ae_a$ for some $d_a\in k[[t]]$ (i.e. independent of $x$).
Again by Lemma \ref{tut}, these projections are also $\tut$-invariant.

It is left to check that the term $z(e_{c_1}+e_{c_2})y$ is central. We show that it commutes with each component of the decomposition \eqref{sumid}.
Indeed, since $e_{c_1}$ and $e_{c_2}$ are $\tut$-invariant (\ref{ei}), then $z(e_{c_1}+e_{c_2})y$ commutes both with
$[{kC_{2^{n-1}}}]_t[y;\tut]e_{c_1}$ and $[{kC_{2^{n-1}}}]_t[y;\tut]e_{c_2}$.
Furthermore, for every $a\in \mathbb{F}^*_{2^{n-2}}$, $e_a$ is orthogonal to $e_{c_1}$ as well as to $e_{c_2}$. Thus,
$$z(e_{c_1}+e_{c_2})y\cdot[{kC_{2^{n-1}}}]_t[y;\tut]e_a=0=[{kC_{2^{n-1}}}]_t[y;\tut]e_a\cdot z(e_{c_1}+e_{c_2})y,$$
and hence $z(e_{c_1}+e_{c_2})y$ commutes with $[{kC_{2^{n-1}}}]_t[y;\tut]$. \qed
%\end{document}
\section{Dihedral 2-groups}\label{dihed}
A slight modification of above construction yields a separable deformation of $kD_{2^n}$, where $D_{2^n}$ is the dihedral group of order $2^n$ (and $\mathbb{F}_{2^{n-2}}\subseteq k$ as before).
It is outlined herein only briefly since the DF conjecture for this family of 2-groups has already been solved in \cite{ES}.

The dihedral group $D_{2^n}$ admits a {\it split} extension
%\begin{equation}\label{dihedextension}
$1\to \Cn\to \Dn\to C_2\to 1,$
%\end{equation}
with the same action $\eta(\btu):\sigma\mapsto \sigma^{2^{n-1}-1}$ of $C_2=\la \bar\tau\ra  $ on $\Cn=\la \sigma\ra  $ as in the extension \eqref{extension}.
Then the dihedral group algebra satisfies $$k\Dn\simeq k\Cn[y;\eta]/\la y^2-1\ra.  $$
We deform the group algebra $k\Cn$ and the action $\eta$ exactly as in \S\ref{C4.1} and \S\ref{C4.2} respectively so as to obtain the deformed skew polynomial algebra
$[{kC_{2^{n-1}}}]_t[y;\tut].$
The difference from the construction in \S\ref{defor} for the generalized quaternions is manifested in the polynomial
$y^2+z(e_{c_1}+e_{c_2})y+1\in[{kC_{2^{n-1}}}]_t[y;\tut],$ where $e_{c_1}$ and $e_{c_2}$ are the primitive idempotents as in \S\ref{C4.1}, and $z\in k[[t]]$ is the same as in \S\ref{hbt}.
This polynomial replaces \eqref{qt1} having a different free term.
Its centrality in $[{kC_{2^{n-1}}}]_t[y;\tut]$ follows from Lemma \ref{gcentproof} (the free term here is obviously central).
Define
$$[k\Dn]_t:=[{kC_{2^{n-1}}}]_t[y;\tut]/\la y^2+z(e_{c_1}+e_{c_2})y+1\ra.$$
Then it is easy to verify that $[k\Dn]_t$ is indeed a deformation of $k\Dn$. Separability of $[k\Dn]_t$ is proven similarly to \S\ref{sep}, moreover we establish again a strong solution
$$\overline{k((t))}\otimes_{k((t))}[kD_{2^{n}}]_t \simeq\bigoplus_{i=1}^{{2^{n-2}-1}}M_2(\overline{k((t))})
\oplus \bigoplus_{i=1}^4\overline{k((t))}.$$

\noindent{\bf Acknowledgement.} The author thanks A. Amsalem for useful discussions.

\end{document}